\documentclass{amsart}
\usepackage{amsmath}
\usepackage{amsfonts}
\usepackage{cite}

\newtheorem{theorem}{Theorem}
\theoremstyle{plain}
\newtheorem{acknowledgement}{Acknowledgement}

\newtheorem{proposition}{Proposition}
\newtheorem{remark}{Remark}

\numberwithin{equation}{section}
\input{tcilatex}

\begin{document}
\title[Killing magnetic curves in non-flat Lorentzian-Heisenberg spaces]{%
Killing magnetic curves in non-flat Lorentzian-Heisenberg spaces}
\author{Murat Altunba\c{s}}
\address[M. Altunbas]{Department of Mathematics \\
Faculty of Sciences and Arts\\
Erzincan University\\
25240 Erzurum, Turkey}
\email[M. Altunbas]{maltunbas@erzincan.edu.tr}

\begin{abstract}
We obtain some explicit formulas for Killing magnetic curves in non-flat
Lorentzian-Heisenberg spaces.

\end{abstract}

\maketitle

\section{\protect\bigskip Introduction}

Let $(M,g)$ be a three-dimensional (semi-)Riemannian manifold.$\ $The
magnetic curves $\gamma $ on $M$ are generalizations of geodesics which
satisfy the following differential equation: 
\begin{equation}
\nabla _{\gamma ^{\prime }}\gamma ^{\prime }=\varphi (\gamma ^{\prime }),
\label{m1}
\end{equation}%
where $\nabla $ is the Levi-Civita connection of $g$ and $\varphi $ is a
skew-symmetric $(1,1)-$tensor field. The tensor field $\varphi $ is known as
the Lorentz force and Equation (\ref{m1}) is said to be the Lorentz
equation. Magnetic curves were investigated by several authors in Riemannian
and semi-Riemannian manifolds (see \cite{Ozdemir}, \cite{Ozgur}, \cite%
{Romaniuc1},\ \cite{Romaniuc2}).

Moreover, when the magnetic fields (which will be explained later)
correspond to a Killing vector, the curves $\gamma $ which fulfill Equation (%
\ref{m1}) are said to be Killing magnetic curves. Studying Killing magnetic
curves is an actual topic of research in pure mathematics and theoretical
physics. In \cite{Romaniuc3}, Romaniuc and Munteanu considered Killing
magnetic curves in three-dimensional Euclidean space. In \cite{Romaniuc4},
the same authors studied these curves in three-dimensional Minkowski space.
In \cite{Erjavec1}, Erjavec gave some characterizations about Killing
magnetic curves in $SL(2,%
\mathbb{R}
).\ $In \cite{Erjavec2} and \cite{Erjavec3}, Killing magnetic curves were
investigated in Sol space and almost cosymplectic Sol space, respectively.
In \cite{Munteanu}, Munteanu and Nistor classified Killing magnetic curves
in $S^{2}\times 
\mathbb{R}
.$ In \cite{Calvaruso}, Calvaruso, Munteanu and Perrone obtained a complete
classification for the Killing magnetic curves in three-dimensional almost
paracontact manifolds. In \cite{Bejan}, Bejan and Romaniuc proved that
equipped with a Killing vector field $V$, any arc length parameterized
spacelike or timelike curve, normal to $V$, is a magnetic trajectory
associated to $V$ in a Walker manifold. And finally, in \cite{Derkaoui},
Derkaoui and Hathout occured explicit formulas for Killing magnetic curves
in Heisenberg group.

In this paper, we determine the Killing magnetic curves in the
three-dimensional Lorentzian-Heisenberg space. It is known that
Lorentzian-Heisenberg space can be equipped with three non-isometric
metrics. We will consider two of them which are non-flat.

\section{Preliminaries}

Let $(M,g)$ be a three-dimensional semi-Riemannian manifold. The magnetic
curves on $(M,g)$ are trajectories of charged particles, moving on $M$ under
the action of a magnetic field $F$. A magnetic field on $M$ is a closed
2-form $F$ on $M$,\ to which one can associate a skew-symmetric $(1,1)-$%
tensor field $\varphi $ on $M$, uniquely determined by 
\begin{equation*}
F(X,Y)=g(\varphi (X),Y),
\end{equation*}%
for all $X,Y\in \chi (M).$ Here, the tensor field $\varphi $ is called the
Lorentz force.

The magnetic trajectories of $F$ are regular curves $\gamma $ in $M$ which
satisfy the Lorentz equation%
\begin{equation}
\nabla _{\mathbf{t}}\mathbf{t}=\varphi (\mathbf{t}),  \label{m2}
\end{equation}%
where $\mathbf{t}=\gamma ^{\prime }\ $is the speed vector of $\gamma .$

Furthermore, to have a positively oriented orthonormal frame field $%
\{e_{1},e_{2},e_{3}\}$ and represent the vectors $X$ and $Y$ as $%
X=x_{1}e_{1}+x_{2}e_{2}+x_{3}e_{3}$\ and $%
Y=y_{1}e_{1}+y_{2}e_{2}+y_{3}e_{3},\ $the vector product of two vector
fields $X=(x_{1},x_{2},x_{3})$ and $Y=(y_{1},y_{2},y_{3})$ is given by 
\begin{equation*}
X\wedge
Y=(x_{2}y_{3}-x_{3}y_{2},x_{3}y_{1}-x_{1}y_{3},x_{2}y_{1}-x_{1}y_{2}).
\end{equation*}%
The mixed product of the vector fields $X,Y,Z\in \chi (M)$ is then defined by%
\begin{equation*}
g(X\wedge Y,Z)=dv_{g}(X,Y,Z),
\end{equation*}%
where $dv_{g\text{ }}$denotes a volume on $M$.

A vector field $V$ is called a Killing vector field if it satisfy the
Killing equation 
\begin{equation*}
g(\nabla _{X}V,Y)+g(\nabla _{Y}V,X)=0,
\end{equation*}%
for all $X,Y\in \chi (M),$ where $\nabla $ is the Levi-Civita connection of
the metric $g.$

Let $F_{V}=i_{V}dvg$ be the Killing magnetic force corresponding to the
Killing magnetic vector field $V$ on $M,$ where $i$ denotes the inner
product. The Lorentz force of $F_{V}\ $is described as 
\begin{equation*}
\varphi (X)=V\wedge X,
\end{equation*}%
for all $X\in \chi (M).\ $Therefore, Equation (\ref{m2}) can be rewritten as 
\begin{equation}
\nabla _{\mathbf{t}}\mathbf{t}=V\wedge \mathbf{t},  \label{m3}
\end{equation}%
and solutions of above equation are called Killing magnetic curves
corresponding to the Killing vector fields $V.$

For shortness, we will call these curves as $V-$magnetic curves in this
paper.

\section{Metrics on Lorentzian-Heisenberg spaces}

Each left-invariant Lorentzian metric on the 3-dimensional Heisenberg group $%
H_{3}$ is isometric to one of the following metrics:%
\begin{eqnarray*}
g_{1} &=&-\frac{1}{\lambda ^{2}}dx^{2}+dy^{2}+(xdy+dz)^{2}, \\
g_{2} &=&\frac{1}{\lambda ^{2}}dx^{2}+dy^{2}-(xdy+dz)^{2},\ \lambda >0, \\
g_{3} &=&dx^{2}+(xdy+dz)^{2}-((1-x)dy-dz)^{2}.
\end{eqnarray*}%
Furthermore, the Lorentzian metrics $g_{1},g_{2},g_{3}$ are non-isometric
and the Lorentzian metric $g_{3}$ is flat \cite{Batat}. We will deal with
the metrics $g_{1}$ and $g_{2}\ ($i.e. non-flat cases).

\begin{remark}
According to the coordinate change $u=\lambda ^{-1}x,\ v=y,\ w=z+2xy,\ $we
rewrite the metrics as%
\begin{eqnarray}
g_{1} &=&-du^{2}+dv^{2}+\lambda ^{2}(vdu-udv)^{2},  \notag \\
g_{2} &=&du^{2}+dv^{2}-\lambda ^{2}(vdu-udv)^{2},\ \lambda >0.  \label{z1}
\end{eqnarray}
\end{remark}

\section{The metric $g_{1}$}

An orthonormal basis on $(H_{3},g_{1})$ is given by 
\begin{equation}
e_{1}=\frac{\partial }{\partial z},\ e_{2}=\frac{\partial }{\partial y}-x%
\frac{\partial }{\partial z},\ e_{3}=\lambda \frac{\partial }{\partial x},
\label{1}
\end{equation}%
where the vector $e_{3}\ $is timelike. The non-zero components of the
Levi-Civita connection $\nabla $ of the metric $g_{1}\ $are given by%
\begin{eqnarray}
\nabla _{e_{1}}e_{2} &=&\nabla _{e_{2}}e_{1}=\frac{\lambda }{2}e_{3},
\label{1a} \\
\nabla _{e_{1}}e_{3} &=&\nabla _{e_{3}}e_{1}=\frac{\lambda }{2}e_{2},  \notag
\\
\nabla _{e_{2}}e_{3} &=&-\nabla _{e_{3}}e_{2}=\frac{\lambda }{2}e_{1}. 
\notag
\end{eqnarray}%
The Lie algebra of Killing vector fields of $(H_{3},g_{1})$ admits as basis 
\begin{eqnarray}
V_{1} &=&\frac{\partial }{\partial z},V_{2}=\frac{\partial }{\partial y},\
V_{3}=\lambda \frac{\partial }{\partial x}-\lambda y\frac{\partial }{%
\partial z},  \label{2} \\
V_{4} &=&\lambda ^{2}y\frac{\partial }{\partial x}+x\frac{\partial }{%
\partial y}-\frac{1}{2}(x^{2}+\lambda ^{2}y^{2})\frac{\partial }{\partial z}.
\notag
\end{eqnarray}%
Using (\ref{1}), we rewrite equations (\ref{2}) as follows:%
\begin{eqnarray*}
V_{1} &=&e_{1},\ V_{2}=xe_{1}+e_{2},V_{3}=-\lambda ye_{1}+e_{3}, \\
V_{4} &=&\frac{1}{2}(x^{2}-\lambda ^{2}y^{2})e_{1}+xe_{2}+\lambda ye_{3}.
\end{eqnarray*}

If $\gamma :I\rightarrow ({\normalsize H}_{3},g_{1}),\ \gamma
(t)=(x(t),y(t),z(t))$ is a regular curve, then its speed vector is described
as 
\begin{equation*}
\mathbf{t}=\gamma ^{\prime }(t)=(x^{\prime }(t),y^{\prime }(t),z^{\prime
}(t))
\end{equation*}%
and%
\begin{equation}
\mathbf{t}=\gamma ^{\prime }(t)=(z^{\prime }+xy^{\prime })e_{1}+y^{\prime
}e_{2}+\frac{x^{\prime }}{\lambda }e_{3}.  \label{c3}
\end{equation}%
From equations (\ref{1a}), we have 
\begin{equation}
\nabla _{\mathbf{t}}\mathbf{t=}(z^{\prime }+xy^{\prime })^{\prime
}e_{1}+(y^{\prime \prime }+x^{\prime }(z^{\prime }+xy^{\prime }))e_{2}+(%
\frac{x^{\prime \prime }}{\lambda }+\lambda y^{\prime }(z^{\prime
}+xy^{\prime }))e_{3}.  \label{3}
\end{equation}%
In the following subsections, we obtain some formulas for $V_{i}\ -$magnetic
curves $(i=1,...,4)$ in $({\normalsize H}_{3},g_{1}).$ To solve the
differential equations, we need help of \textit{Wolfram Mathematica.}

\subsection{$V_{1}-$magnetic curves}

\bigskip Using (\ref{2})$_{1\text{ }}$and (\ref{c3}), we occur 
\begin{equation}
V_{1}\wedge \mathbf{t=-}\frac{x^{\prime }}{\lambda }e_{2}-y^{\prime }e_{3}.
\label{4}
\end{equation}

The equation $\nabla _{\mathbf{t}}\mathbf{t=}V_{1}\wedge \mathbf{t}$ gives
us the following system: 
\begin{equation}
S_{1}:\left\{ 
\begin{array}{l}
y^{\prime \prime }+x^{\prime }((z^{\prime }+xy^{\prime })+\frac{1}{\lambda }%
)=0, \\ 
x^{\prime \prime }+y^{\prime }(\lambda ^{2}(z^{\prime }+xy^{\prime
})+\lambda )=0, \\ 
(z^{\prime }+xy^{\prime })^{\prime }=0.%
\end{array}%
\right.  \label{5}
\end{equation}%
By integrating $(S_{1})_{3}$\ and putting it in $(S_{1})_{1,2},$ we obtain 
\begin{equation*}
S_{1}:\left\{ 
\begin{array}{l}
y^{\prime \prime }+x^{\prime }(c+\frac{1}{\lambda })=0, \\ 
x^{\prime \prime }+y^{\prime }\lambda (\lambda c+1)=0, \\ 
z^{\prime }+xy^{\prime }=c\ \text{(constant).}%
\end{array}%
\right.
\end{equation*}%
Solution of the system $(S_{1})_{1,2}$ is 
\begin{equation}
\left\{ 
\begin{array}{l}
x(t)=-\frac{\lambda }{\lambda c+1}[k_{1}\cosh ((\lambda c+1)t)+k_{2}\sinh
((\lambda c+1)t)]+k_{3}, \\ 
y(t)=\frac{1}{\lambda c+1}[k_{1}\sinh ((\lambda c+1)t)+k_{2}\cosh ((\lambda
c+1)t)]+k_{4},%
\end{array}%
\right.  \label{6}
\end{equation}%
where $k_{i},\ i=1,...,4\ $are constants. Setting equations (\ref{6}) in $%
(S_{1})_{3}$ and by integration, we get%
\begin{eqnarray*}
z(t) &=&(c+\frac{\lambda }{2(\lambda c+1)}(k_{1}^{2}-k_{2}^{2}))t+\frac{%
\lambda }{(\lambda c+1)^{2}}[\frac{(k_{1}^{2}+k_{2}^{2})}{4}\sinh (2(\lambda
c+1)t) \\
&&+\frac{k_{1}k_{2}}{2}\cosh (2(\lambda c+1)t)]+\frac{1}{\lambda c+1}%
[k_{1}k_{3}\sinh ((\lambda c+1)t) \\
&&+k_{2}k_{3}\cosh ((\lambda c+1)t)]+k_{5},
\end{eqnarray*}%
where $k_{5\text{ }}$is a constant.

If $c=-\frac{1}{\lambda },\ $the system $S_{1}$ reduces to 
\begin{equation*}
S_{1}:\left\{ 
\begin{array}{l}
y^{\prime \prime }=0, \\ 
x^{\prime \prime }=0, \\ 
z^{\prime }+xy^{\prime }=-\frac{1}{\lambda }\text{.}%
\end{array}%
\right.
\end{equation*}%
Its general solution 
\begin{equation*}
S_{1}:\left\{ 
\begin{array}{l}
x(t)=k_{1}t+k_{2}, \\ 
y(t)=k_{3}t+k_{4}, \\ 
z(t)=-\frac{k_{1}k_{3}}{2}t^{2}-(\frac{1}{\lambda }+k_{2}k_{3})t+k_{5}\text{,%
}%
\end{array}%
\right.
\end{equation*}%
where $k_{i},\ i=1,...,5\ $are constants. Therefore, we state the following
theorem.

\begin{theorem}
All $V_{1}-$magnetic curves of $({\normalsize H}_{3},g_{1})$ satisfy the
following equations:

(i)\ If$\ c=-\frac{1}{\lambda },\ $then 
\begin{equation*}
\gamma (t)=\left( 
\begin{array}{c}
x(t)=k_{1}t+k_{2}, \\ 
y(t)=k_{3}t+k_{4}, \\ 
z(t)=-\frac{k_{1}k_{3}}{2}t^{2}-(\frac{1}{\lambda }+k_{2}k_{3})t+k_{5}%
\end{array}%
\right) .
\end{equation*}

(ii) If$\ c\neq -\frac{1}{\lambda },\ $then%
\begin{equation*}
\gamma (t)=\left( 
\begin{array}{c}
x(t)=-\frac{\lambda }{\lambda c+1}[k_{1}\cosh ((\lambda c+1)t)+k_{2}\sinh
((\lambda c+1)t)]+k_{3}, \\ 
y(t)=\frac{1}{\lambda c+1}[k_{1}\sinh ((\lambda c+1)t)+k_{2}\cosh ((\lambda
c+1)t)]+k_{4}, \\ 
z(t)=(c+\frac{\lambda }{2(\lambda c+1)}(k_{1}^{2}-k_{2}^{2}))t+\frac{\lambda 
}{(\lambda c+1)^{2}}[\frac{(k_{1}^{2}+k_{2}^{2})}{4}\sinh (2(\lambda c+1)t)
\\ 
+\frac{k_{1}k_{2}}{2}\cosh (2(\lambda c+1)t)]+\frac{1}{(\lambda c+1)}%
[k_{1}k_{3}\sinh ((\lambda c+1)t) \\ 
+k_{2}k_{3}\cosh ((\lambda c+1)t)]+k_{5}%
\end{array}%
\right) ,
\end{equation*}%
where $k_{i},\ i=1,...,5\ $are constants.
\end{theorem}

\subsection{$V_{2}-$magnetic curves}

According to (\ref{2})$_{2\text{ }}$and (\ref{c3}), we have%
\begin{equation}
V_{2}\wedge \mathbf{t=}\frac{x^{\prime }}{\lambda }e_{1}-\frac{xx^{\prime }}{%
\lambda }e_{2}-(xy^{\prime }-(z^{\prime }+xy^{\prime }))e_{3}.  \label{7}
\end{equation}%
From the equation $\nabla _{\mathbf{t}}\mathbf{t=}V_{2}\wedge \mathbf{t,\ }$%
we get%
\begin{equation}
S_{2}:\left\{ 
\begin{array}{l}
y^{\prime \prime }+x^{\prime }(z^{\prime }+xy^{\prime })=-\frac{xx^{\prime }%
}{\lambda }, \\ 
(\lambda y^{\prime }-1)(z^{\prime }+xy^{\prime })=-xy^{\prime }-\frac{%
x^{\prime \prime }}{\lambda }, \\ 
(z^{\prime }+xy^{\prime })^{\prime }=\frac{x^{\prime }}{\lambda }.%
\end{array}%
\right.  \label{8}
\end{equation}

By integrating $(S_{2})_{3},\ $we deduce%
\begin{equation*}
z^{\prime }+xy^{\prime }=\frac{x}{\lambda }+c,
\end{equation*}%
where $c$ is a constant. Putting the last equation in $(S_{2})_{1,2}$ we get%
\begin{equation*}
\bar{S}_{2}:\left\{ 
\begin{array}{l}
y^{\prime \prime }=-x^{\prime }(\frac{2x}{\lambda }+c), \\ 
(\lambda y^{\prime }-1)(\frac{x}{\lambda }+c)=-xy^{\prime }-\frac{x^{\prime
\prime }}{\lambda },%
\end{array}%
\right.
\end{equation*}%
and%
\begin{equation}
\bar{S}_{2}:\left\{ 
\begin{array}{l}
y^{\prime }=-\frac{x^{2}}{\lambda }-xc, \\ 
x^{\prime \prime }-2x^{3}-x-\lambda c(3x^{2}+c\lambda x+1)=0.%
\end{array}%
\right.  \label{9}
\end{equation}%
Without loss of generality, we can suppose that $c=0.$\ In this case\emph{,}$%
\ $the equation $x^{\prime \prime }-2x^{3}-x=0$\emph{\ }involves Jacobi
elliptic functions as solutions. So, we can express the following
proposition.

\begin{proposition}
The Killing magnetic curves in $(H_{3},g_{1})$ corresponding to the Killing
vector field $V_{2}=xe_{1}+e_{2}$ are solutions of the system of
differential equations (\ref{8}).
\end{proposition}

\subsection{$V_{3}-$magnetic curves}

From (\ref{2})$_{3\text{ }}$and (\ref{c3}), we have 
\begin{equation}
V_{3}\wedge \mathbf{t=-}y^{\prime }e_{1}+(yx^{\prime }+(z^{\prime
}+xy^{\prime }))e_{2}+\lambda yy^{\prime }e_{3}.  \label{10}
\end{equation}%
Using the equation $\nabla _{\mathbf{t}}\mathbf{t=}V_{3}\wedge \mathbf{t,\ }$%
we obtain%
\begin{equation}
S_{3}:\left\{ 
\begin{array}{l}
y^{\prime \prime }+x^{\prime }(z^{\prime }+xy^{\prime })=yx^{\prime
}+(z^{\prime }+xy^{\prime }), \\ 
\frac{x^{\prime \prime }}{\lambda }+\lambda y^{\prime }(z^{\prime
}+xy^{\prime })=\lambda yy^{\prime }, \\ 
(z^{\prime }+xy^{\prime })^{\prime }=-y^{\prime }.%
\end{array}%
\right.  \label{11}
\end{equation}

By integrating $(S_{3})_{3},\ $we occur%
\begin{equation*}
z^{\prime }+xy^{\prime }=-y+c,
\end{equation*}%
where $c$ is a constant. Putting the last equation in $(S_{3})_{1,2},$ we get%
\begin{equation}
\bar{S}_{3}:\left\{ 
\begin{array}{l}
x^{\prime }=\lambda ^{2}y^{2}-\lambda ^{2}yc, \\ 
y^{\prime \prime }-2x^{\prime }y+y+c(x^{\prime }-1)=0.%
\end{array}%
\right.  \label{12}
\end{equation}%
Without loss of generality, we can assume that $c=0.\ $In this case\emph{,}$%
\ $when we try to solve the system $\bar{S}_{3},\ i.e.,\ $the equation\ $%
y^{\prime \prime }-2\lambda ^{2}y^{3}+y=0,$\emph{\ }we encounter Jacobi
elliptic functions. Therefore, we write the following proposition.

\begin{proposition}
The Killing magnetic curves in $(H_{3},g_{1})$ corresponding to the Killing
vector field $V_{3}=-\lambda ye_{1}+e_{3}$ are solutions of the system of
differential equations (\ref{11}).
\end{proposition}

\subsection{$V_{4}-$magnetic curves}

From (\ref{2})$_{3\text{ }}$and (\ref{c3}), we write 
\begin{eqnarray}
V_{4}\wedge \mathbf{t} &=&(\frac{xx^{\prime }}{\lambda }-\lambda yy^{\prime
})e_{1}-(\frac{1}{2}(x^{2}-\lambda ^{2}y^{2})\frac{x^{\prime }}{\lambda }%
-\lambda y(z^{\prime }+xy^{\prime }))e_{2}  \label{13} \\
&&-(\frac{1}{2}(x^{2}-\lambda ^{2}y^{2})y^{\prime }-x(z^{\prime }+xy^{\prime
}))e_{3}.  \notag
\end{eqnarray}%
From the equation $\nabla _{\mathbf{t}}\mathbf{t=}V_{4}\wedge \mathbf{t,\ }$%
we get%
\begin{equation}
S_{4}:\left\{ 
\begin{array}{l}
y^{\prime \prime }+x^{\prime }(z^{\prime }+xy^{\prime })=-\frac{1}{2}%
(x^{2}-\lambda ^{2}y^{2})\frac{x^{\prime }}{\lambda }+\lambda y(z^{\prime
}+xy^{\prime }), \\ 
\frac{x^{\prime \prime }}{\lambda }+\lambda y^{\prime }(z^{\prime
}+xy^{\prime })=-\frac{1}{2}(x^{2}-\lambda ^{2}y^{2})y^{\prime }+x(z^{\prime
}+xy^{\prime }), \\ 
(z^{\prime }+xy^{\prime })^{\prime }=(\frac{xx^{\prime }}{\lambda }-\lambda
yy^{\prime }).%
\end{array}%
\right.  \label{14}
\end{equation}

By integrating $(S_{4})_{3},\ $we obtain%
\begin{equation}
z^{\prime }+xy^{\prime }=\frac{x^{2}}{2\lambda }-\frac{\lambda y^{2}}{2}+c,
\label{144}
\end{equation}%
where $c$ is a constant. Putting the last equation in $(S_{4})_{1,2},$ we get%
\begin{equation}
\bar{S}_{4}:\left\{ 
\begin{array}{l}
y^{\prime \prime }+\frac{x^{\prime }}{\lambda }(x^{2}-\lambda ^{2}y^{2})=%
\frac{1}{2}y(x^{2}-\lambda ^{2}y^{2})+c(\lambda y-x^{\prime }), \\ 
x^{\prime \prime }+\lambda y^{\prime }(x^{2}-\lambda ^{2}y^{2})=\frac{1}{2}%
x(x^{2}-\lambda ^{2}y^{2})+c(\lambda x-\lambda ^{2}y^{\prime }).%
\end{array}%
\right.  \label{15}
\end{equation}%
It seems very difficult to solve the system $\bar{S}_{4}$ in general case.
For a particular case\emph{\ }$x=\lambda y$\emph{, }we deduce 
\begin{equation}
\bar{S}_{4}:\left\{ 
\begin{array}{l}
x^{\prime \prime }+c\lambda x^{\prime }-c\lambda x=0, \\ 
y^{\prime \prime }+c\lambda y^{\prime }-c\lambda y=0.%
\end{array}%
\right.  \label{155}
\end{equation}%
By solving the second equation of the above system, we get 
\begin{equation*}
y(t)=k_{1}e^{-\frac{t}{2}(c\lambda +\sqrt{c\lambda (4+c\lambda )})}+k_{2}e^{%
\frac{t}{2}(-c\lambda +\sqrt{c\lambda (4+c\lambda )})},
\end{equation*}%
where $k_{1}$ and $k_{2}$ are constants. From (\ref{144}), we obtain 
\begin{eqnarray*}
z(t) &=&ct-\frac{\lambda y^{2}}{2} \\
&=&ct-\frac{\lambda }{2}(k_{1}e^{-\frac{t}{2}(c\lambda +\sqrt{c\lambda
(4+c\lambda )})}+k_{2}e^{\frac{t}{2}(-c\lambda +\sqrt{c\lambda (4+c\lambda )}%
)})^{2}.
\end{eqnarray*}%
Therefore, the solution of the system $\bar{S}_{4}$ is given by 
\begin{equation}
\bar{S}_{4}:\left\{ 
\begin{array}{l}
x(t)=\lambda \left( k_{1}e^{-\frac{t}{2}(c\lambda +\sqrt{c\lambda
(4+c\lambda )})}+k_{2}e^{\frac{t}{2}(-c\lambda +\sqrt{c\lambda (4+c\lambda )}%
)}\right) , \\ 
y(t)=k_{1}e^{-\frac{t}{2}(c\lambda +\sqrt{c\lambda (4+c\lambda )})}+k_{2}e^{%
\frac{t}{2}(-c\lambda +\sqrt{c\lambda (4+c\lambda )})}, \\ 
z(t)=ct-\frac{\lambda }{2}(k_{1}e^{-\frac{t}{2}(c\lambda +\sqrt{c\lambda
(4+c\lambda )})}+k_{2}e^{\frac{t}{2}(-c\lambda +\sqrt{c\lambda (4+c\lambda )}%
)})^{2}.%
\end{array}%
\right.  \label{156}
\end{equation}%
Hence, we write the following proposition.

\begin{proposition}
The Killing magnetic curves in $(H_{3},g_{1})$ corresponding to the Killing
vector field $V_{4}=\frac{1}{2}(x^{2}-\lambda ^{2}y^{2})e_{1}+xe_{2}+\lambda
ye_{3}$ are solutions of the system of differential equations (\ref{14}).
Moreover, the space curves given by parametric equations (\ref{156}) are $%
V_{4}-$magnetic curves in $(H_{3},g_{1}).$
\end{proposition}

In the last section, we follow the steps explained in the strategy mentioned
in this section for the metric $g_{2}.$

\section{The metric $g_{2}$}

We have an orthonormal basis on $(H_{3},g_{2})$ 
\begin{equation}
\ e_{1}=\frac{\partial }{\partial y}-x\frac{\partial }{\partial z},\
e_{2}=\lambda \frac{\partial }{\partial x},e_{3}=\frac{\partial }{\partial z}%
,  \label{16}
\end{equation}%
where the vector $e_{3}\ $is timelike. The non-zero components of the
Levi-Civita connection $\nabla $ of the metric $g_{2}\ $are given by%
\begin{eqnarray}
\nabla _{e_{1}}e_{2} &=&-\nabla _{e_{2}}e_{1}=\frac{\lambda }{2}e_{3},
\label{166} \\
\nabla _{e_{1}}e_{3} &=&\nabla _{e_{3}}e_{1}=\frac{\lambda }{2}e_{2},  \notag
\\
\nabla _{e_{2}}e_{3} &=&\nabla _{e_{3}}e_{2}=-\frac{\lambda }{2}e_{1}. 
\notag
\end{eqnarray}%
The Lie algebra of Killing vector fields of $(H_{3},g_{2})$ admits as basis 
\begin{eqnarray}
\ V_{1} &=&\frac{\partial }{\partial z},V_{2}=\frac{\partial }{\partial y}%
,V_{3}=\lambda \frac{\partial }{\partial x}-\lambda y\frac{\partial }{%
\partial z},\   \label{17} \\
V_{4} &=&-\lambda ^{2}y\frac{\partial }{\partial x}+x\frac{\partial }{%
\partial y}+\frac{1}{2}(-x^{2}+\lambda ^{2}y^{2})\frac{\partial }{\partial z}%
.  \notag
\end{eqnarray}%
Using (\ref{16}), we rewrite Equations (\ref{17}) as follows:%
\begin{eqnarray*}
V_{1} &=&e_{3},\ V_{2}=e_{1}+xe_{3},V_{3}=e_{2}-\lambda ye_{3}, \\
V_{4} &=&xe_{1}-\lambda ye_{2}+\frac{1}{2}(x^{2}+\lambda ^{2}y^{2})e_{3}.
\end{eqnarray*}

If $\gamma :I\rightarrow ({\normalsize H}_{3},g_{2}),\ \gamma
(t)=(x(t),y(t),z(t))$ is a regular curve, then its speed vector is described
as 
\begin{equation*}
\mathbf{t}=\gamma ^{\prime }(t)=(x^{\prime }(t),y^{\prime }(t),z^{\prime
}(t))
\end{equation*}%
and%
\begin{equation}
\mathbf{t}=\gamma ^{\prime }(t)=y^{\prime }e_{1}+\frac{x^{\prime }}{\lambda }%
e_{2}+(z^{\prime }+xy^{\prime })e_{3}.  \label{m17}
\end{equation}%
From equations (\ref{166}), we have%
\begin{equation}
\nabla _{\mathbf{t}}\mathbf{t=}(y^{\prime \prime }-x^{\prime }(z^{\prime
}+xy^{\prime }))e_{1}+(\frac{x^{\prime \prime }}{\lambda }+\lambda y^{\prime
}(z^{\prime }+xy^{\prime }))e_{2}+(z^{\prime }+xy^{\prime })^{\prime }e_{3}.
\label{18}
\end{equation}

\subsection{$V_{1}-$magnetic curves}

\bigskip We have 
\begin{equation}
V_{1}\wedge \mathbf{t=-}\frac{x^{\prime }}{\lambda }e_{1}+y^{\prime }e_{2}.
\label{19}
\end{equation}

From the equation $\nabla _{\mathbf{t}}\mathbf{t=}V_{1}\wedge \mathbf{t,\ }$%
we get%
\begin{equation}
S_{1}:\left\{ 
\begin{array}{l}
y^{\prime \prime }-x^{\prime }((z^{\prime }+xy^{\prime })-\frac{1}{\lambda }%
)=0, \\ 
x^{\prime \prime }+y^{\prime }(\lambda ^{2}(z^{\prime }+xy^{\prime
})-\lambda )=0, \\ 
(z^{\prime }+xy^{\prime })^{\prime }=0.%
\end{array}%
\right.  \label{20}
\end{equation}%
By integrating $(S_{1})_{3}$\ and putting it in $(S_{1})_{1,2},$ we obtain 
\begin{equation*}
S_{1}:\left\{ 
\begin{array}{l}
y^{\prime \prime }-x^{\prime }(c-\frac{1}{\lambda })=0, \\ 
x^{\prime \prime }+y^{\prime }\lambda (\lambda c-1)=0, \\ 
z^{\prime }+xy^{\prime }=c\ \text{(constant).}%
\end{array}%
\right.
\end{equation*}%
Solution of the system $(S_{1})_{1,2}$ is 
\begin{equation}
S_{1}:\left\{ 
\begin{array}{l}
x(t)=\frac{\lambda }{\lambda c-1}[k_{1}\cos ((\lambda c-1)t)+k_{2}\sin
((\lambda c-1)t)]+k_{3}, \\ 
y(t)=\frac{1}{\lambda c-1}[k_{1}\sin ((\lambda c-1)t)-k_{2}\cos ((\lambda
c-1)t)]+k_{4},%
\end{array}%
\right.  \label{21}
\end{equation}%
where $k_{i},\ i=1,...,4\ $are constants. Setting equations (\ref{20}) in $%
(S_{1})_{3}$ and by integration, we get%
\begin{eqnarray*}
z(t) &=&(c-\frac{\lambda }{2(\lambda c-1)}(k_{1}^{2}-k_{2}^{2}))t-\frac{%
\lambda }{(\lambda c-1)^{2}}[\frac{(k_{1}^{2}-k_{2}^{2})}{4}\sin (2(\lambda
c-1)t) \\
&&-\frac{k_{1}k_{2}}{2}\cos (2(\lambda c-1)t)]+\frac{1}{\lambda c-1}%
(k_{1}k_{3}\sin ((\lambda c-1)t) \\
&&-k_{2}k_{3}\cos ((\lambda c-1)t))+k_{5},
\end{eqnarray*}%
where $k_{5\text{ }}$is a constant. If $c=\frac{1}{\lambda },\ $the system $%
S_{1}$ reduces to 
\begin{equation*}
S_{1}:\left\{ 
\begin{array}{l}
y^{\prime \prime }=0, \\ 
x^{\prime \prime }=0, \\ 
z^{\prime }+xy^{\prime }=\frac{1}{\lambda }\text{.}%
\end{array}%
\right.
\end{equation*}%
Its general solution 
\begin{equation*}
S_{1}:\left\{ 
\begin{array}{l}
x(t)=k_{1}t+k_{2}, \\ 
y(t)=k_{3}t+k_{4}, \\ 
z(t)=-\frac{k_{1}k_{3}}{2}t^{2}+(\frac{1}{\lambda }-k_{2}k_{3})t+k_{5}\text{,%
}%
\end{array}%
\right.
\end{equation*}%
where $k_{i},\ i=1,...,5\ $are constants. So, we have proved the theorem
below.

\begin{theorem}
All $V_{1}-$magnetic curves of $({\normalsize H}_{3},g_{2})$ satisfy the
following equations:

(i)\ If$\ c=\frac{1}{\lambda },\ $then 
\begin{equation*}
\gamma (t)=\left( 
\begin{array}{c}
x(t)=k_{1}t+k_{2}, \\ 
y(t)=k_{3}t+k_{4}, \\ 
z(t)=-\frac{k_{1}k_{3}}{2}t^{2}+(\frac{1}{\lambda }-k_{2}k_{3})t+k_{5}%
\end{array}%
\right) .
\end{equation*}

(ii) If$\ c\neq \frac{1}{\lambda },\ $then%
\begin{equation*}
\gamma (t)=\left( 
\begin{array}{c}
x(t)=\frac{\lambda }{\lambda c-1}[k_{1}\cos ((\lambda c-1)t)+k_{2}\sin
((\lambda c-1)t)]+k_{3}, \\ 
y(t)=\frac{1}{\lambda c-1}[k_{1}\sin ((\lambda c-1)t)-k_{2}\cos ((\lambda
c-1)t)]+k_{4}, \\ 
z(t)=(c-\frac{\lambda }{2(\lambda c-1)}(k_{1}^{2}-k_{2}^{2}))t-\frac{\lambda 
}{(\lambda c-1)^{2}}[\frac{(k_{1}^{2}-k_{2}^{2})}{4}\sin (2(\lambda c-1)t)
\\ 
-\frac{k_{1}k_{2}}{2}\cos (2(\lambda c-1)t)]+\frac{1}{\lambda c-1}%
(k_{1}k_{3}\sin ((\lambda c-1)t) \\ 
-k_{2}k_{3}\cos ((\lambda c-1)t))+k_{5}%
\end{array}%
\right) ,
\end{equation*}%
where $k_{i},\ i=1,...,5\ $are constants.
\end{theorem}

\begin{remark}
These curves were considered by Lee in \cite{Lee} according to corresponding
metric $g_{2}\ $in (\ref{z1}).
\end{remark}

\subsection{$V_{2}-$magnetic curves}

Direct computations give 
\begin{equation}
V_{2}\wedge \mathbf{t=}-\frac{xx^{\prime }}{\lambda }e_{1}+(xy^{\prime
}-(z^{\prime }+xy^{\prime }))e_{2}\mathbf{-}\frac{x^{\prime }}{\lambda }%
e_{3}.  \label{22}
\end{equation}%
The equation $\nabla _{\mathbf{t}}\mathbf{t=}V_{2}\wedge \mathbf{t}$
concludes%
\begin{equation}
S_{2}:\left\{ 
\begin{array}{l}
y^{\prime \prime }-x^{\prime }(z^{\prime }+xy^{\prime })=-\frac{xx^{\prime }%
}{\lambda }, \\ 
(\lambda y^{\prime }+1)(z^{\prime }+xy^{\prime })=xy^{\prime }-\frac{%
x^{\prime \prime }}{\lambda }, \\ 
(z^{\prime }+xy^{\prime })^{\prime }=-\frac{x^{\prime }}{\lambda }.%
\end{array}%
\right.  \label{23}
\end{equation}

By integrating $(S_{2})_{3},\ $we obtain%
\begin{equation*}
z^{\prime }+xy^{\prime }=-\frac{x}{\lambda }+c,
\end{equation*}%
where $c$ is a constant. Putting the last equation in $(S_{2})_{1,2},$ we get%
\begin{equation}
\bar{S}_{2}:\left\{ 
\begin{array}{l}
y^{\prime }=-\frac{x^{2}}{\lambda }+xc, \\ 
x^{\prime \prime }+2x^{3}-x-\lambda c(3x^{2}-c\lambda x-1)=0.%
\end{array}%
\right.  \label{24}
\end{equation}%
Without loss of generality, we can assume that $c=0.$\ This system$\ \bar{S}%
_{2},\ i.e.,\ $the equation $x^{\prime \prime }+2x^{3}-x=0$\emph{\ }involves
Jacobi elliptic functions. So, we write the following proposition.

\begin{proposition}
The Killing magnetic curves in $(H_{3},g_{2})$ corresponding to the Killing
vector field $V_{2}=e_{1}+xe_{3}$ are solutions of the system of
differential equations (\ref{23}).
\end{proposition}

\subsection{\protect\bigskip $V_{3}-$magnetic curves}

We have 
\begin{equation}
V_{3}\wedge \mathbf{t=}(yx^{\prime }+(z^{\prime }+xy^{\prime
}))e_{1}-\lambda yy^{\prime }e_{2}+y^{\prime }e_{3}.  \label{25}
\end{equation}%
From the equation $\nabla _{\mathbf{t}}\mathbf{t=}V_{3}\wedge \mathbf{t,\ }$%
we get%
\begin{equation}
S_{3}:\left\{ 
\begin{array}{l}
y^{\prime \prime }-x^{\prime }(z^{\prime }+xy^{\prime })=yx^{\prime
}+(z^{\prime }+xy^{\prime }), \\ 
\frac{x^{\prime \prime }}{\lambda }+\lambda y^{\prime }(z^{\prime
}+xy^{\prime })=-\lambda yy^{\prime }, \\ 
(z^{\prime }+xy^{\prime })^{\prime }=y^{\prime }.%
\end{array}%
\right.  \label{26}
\end{equation}

By integrating $(S_{3})_{3},\ $we deduce%
\begin{equation*}
z^{\prime }+xy^{\prime }=y+c,
\end{equation*}%
where $c$ is a constant. Putting the last equation in $(S_{3})_{1,2},$ we
have%
\begin{equation}
\bar{S}_{3}:\left\{ 
\begin{array}{l}
x^{\prime }=-\lambda ^{2}y^{2}-\lambda ^{2}yc, \\ 
y^{\prime \prime }-2x^{\prime }y-y-c(x^{\prime }+1)=0.%
\end{array}%
\right.  \label{27}
\end{equation}%
Without loss of generality, we can suppose that $c=0.$\ Then\emph{, }the
system $\bar{S}_{3},\ i.e.,\ $the equation $\ y^{\prime \prime }+2\lambda
^{2}y^{3}-y=0$\emph{\ }has solutions which include Jacobi elliptic functions%
\emph{. }Thus, we give the proposition below.

\begin{proposition}
The Killing magnetic curves in $(H_{3},g_{2})$ corresponding to the Killing
vector field $V_{3}=e_{2}-\lambda ye_{3}$ are solutions of the system of
differential equations (\ref{26}).
\end{proposition}

\subsection{$V_{4}-$magnetic curves}

We write 
\begin{equation}
V_{4}\wedge \mathbf{t=(}-\frac{1}{2}(x^{2}+\lambda ^{2}y^{2})\frac{x^{\prime
}}{\lambda }-\lambda y(z^{\prime }+xy^{\prime }))e_{1}+((x^{2}+\lambda
^{2}y^{2})\frac{y^{\prime }}{2}-x(z^{\prime }+xy^{\prime }))e_{2}-(\frac{%
xx^{\prime }}{\lambda }+\lambda yy^{\prime })e_{3}.  \label{28}
\end{equation}%
The equation $\nabla _{\mathbf{t}}\mathbf{t=}V_{4}\wedge \mathbf{t}$ gives us%
\begin{equation}
S_{4}:\left\{ 
\begin{array}{l}
y^{\prime \prime }-x^{\prime }(z^{\prime }+xy^{\prime })=-\frac{1}{2}%
(x^{2}+\lambda ^{2}y^{2})\frac{x^{\prime }}{\lambda }-\lambda y(z^{\prime
}+xy^{\prime }), \\ 
\frac{x^{\prime \prime }}{\lambda }+\lambda y^{\prime }(z^{\prime
}+xy^{\prime })=\frac{1}{2}(x^{2}+\lambda ^{2}y^{2})y^{\prime }-x(z^{\prime
}+xy^{\prime }), \\ 
(z^{\prime }+xy^{\prime })^{\prime }=-\frac{xx^{\prime }}{\lambda }-\lambda
yy^{\prime }.%
\end{array}%
\right.  \label{29}
\end{equation}

By integrating $(S_{4})_{3},\ $we obtain%
\begin{equation}
z^{\prime }+xy^{\prime }=-\frac{x^{2}}{2\lambda }-\frac{\lambda y^{2}}{2}+c,
\label{30}
\end{equation}%
where $c$ is a constant. Putting the last equation in $(S_{4})_{1,2},$ we get%
\begin{equation}
\bar{S}_{4}:\left\{ 
\begin{array}{l}
y^{\prime \prime }+\frac{x^{\prime }}{\lambda }(x^{2}+\lambda ^{2}y^{2})=%
\frac{1}{2}y(x^{2}+\lambda ^{2}y^{2})+c(-\lambda y+x^{\prime }), \\ 
x^{\prime \prime }-\lambda y^{\prime }(x^{2}+\lambda ^{2}y^{2})=\frac{1}{2}%
x(x^{2}+\lambda ^{2}y^{2})-c(\lambda x+\lambda ^{2}y^{\prime }).%
\end{array}%
\right.  \label{31}
\end{equation}%
It seems a true challenge to solve the system $\bar{S}_{4}$ in general case.
However, we can find a special solution by considering $c=\lambda =1.$ In
this case,%
\begin{equation*}
x(t)=\cos \frac{\sqrt{2}}{2}t,\ y(t)=\sin \frac{\sqrt{2}}{2}t
\end{equation*}%
will be a solution for the system $\bar{S}_{4}.\ $Using these relations in (%
\ref{30}), we get 
\begin{equation*}
z(t)=\frac{2-\sqrt{2}}{4}t-\frac{1}{4}\sin \sqrt{2}t+k_{1}
\end{equation*}%
where $k_{1}$ is a constant. So, we find a solution as follows: 
\begin{equation}
\left\{ 
\begin{array}{l}
x(t)=\cos \frac{\sqrt{2}}{2}t, \\ 
\ y(t)=\sin \frac{\sqrt{2}}{2}t, \\ 
z(t)=\frac{2-\sqrt{2}}{4}t-\frac{1}{4}\sin \sqrt{2}t+k_{1}.%
\end{array}%
\right.  \label{32}
\end{equation}%
Therefore, we can state the last propositon of the paper.

\begin{proposition}
The space curves given by parametric equations%
\begin{equation*}
\gamma (t)=\left( 
\begin{array}{c}
x(t)=\cos \frac{\sqrt{2}}{2}t, \\ 
\ y(t)=\sin \frac{\sqrt{2}}{2}t, \\ 
z(t)=\frac{2-\sqrt{2}}{4}t-\frac{1}{4}\sin \sqrt{2}t+k_{1}%
\end{array}%
\right)
\end{equation*}%
are $V_{4}-$magnetic curves in $(H_{3},g_{2}),\ $where $k_{1}$ is an
arbitrary constants.
\end{proposition}

\begin{acknowledgement}
The author would like to thank Professor Zlatko Erjavec and Professor
Jun-ichi Inoguchi for their valuable suggestions.
\end{acknowledgement}

\end{document}